# A Distance Measuring Algorithm for Location Analysis

**Ruilin Ouyang, Dinghao Ma, M. S. Morshed and Md. Noor-E-Alam\***

Department of Mechanical and Industrial Engineering

Northeastern University

360 Huntington Ave, Boston, MA 02115, USA

**Abstract:**

Approximating distance is one of the key challenge in a facility location problem. Several algorithms have been proposed, however, none of them focused on estimating distance between two concave regions. In this work, we present an algorithm to estimate the distance between two irregular regions of a facility location problem. The proposed algorithm can identify the distance between concave shape regions. We also discuss some relevant properties of the proposed algorithm. A distance-sensitive capacity location model is introduced to test the algorithm. Moreover, several special geometric cases are discussed to show the advantages and insights of the algorithm.

Key words: Distance measurement, facility location problem, $l_p$ distance, geometric distance.

## 1. Introduction:

Calculating distance between two irregular shape regions is known as a key challenge in the context of facility location problem. It's very easy to calculate the geometrical distance between two points. However, when we consider the distance between two irregular regions, it's difficult to say that we have an accurate solution. Several algorithms have been developed to estimate the distance between regions - some of them only works well under a specific condition.

*P*-norm distance or "$l_p$ distance" is widely used in the facility location problem to calculate distance between two points. The following equation for "$l_p$ distance" was proposed by Francis (1967).

\*Corresponding author email: mnalam@neu.edu

$$l_p(q, r) = \left[\sum_{i=1}^{n} |q_i - r_i|^p\right]^{1/p} = \|q - r\|_p = \|\vec{rq}\|_p$$

; where $q = (q_1, \ldots, q_n)$ and $r = (r_1, \ldots, r_n)$ are two points in $n$-dimensional space (i.e. $q, r \in \mathbb{R}^n$). Love and Morris (1972) compared seven different distance functions including $l_p$ functions, weighted $l_p$ functions and weighted Euclidean functions on two road distance samples to determine the accuracy of various distance functions. They considered two criterions with respect to the accuracy of estimation, which are:

- Sum of absolute deviations (the difference between actual road miles and estimated road miles)
- Sum of squares (more sensitive and require more accuracy on shorter distance).

Love and Morris (1979) further compared different parameters for the "$l_p$ distance" model to make the distance calculation more accurate. They used statistical tests to compare the results on seven samples and concluded that "the more parameters used in the model, the more accurate and complexity will get". Love and Dowling (1985) used "$l_{k,p}$ distance" model for floor layout problems to find the appropriate *k, p* values compared to the "$l_1$ distance" (or rectangular distance) model. In the literature, authors stated that the "$l_1$ distance" model was suitable for obtaining optimal facility locations while "$l_{k,p}$ distance" model should be used when total cost function is considered by adjusting *k* values rather than *p* values. Despite several generalizations of this algorithm, none of these algorithms mentioned can successfully determine the exact location of the potential point is in each region.

Batta et al. (1989) solved the location problem with convex forbidden regions using Manhattan metric *(p=1)*. Amiri-Aref et al. (2011) used *P*-norm distance model to find the new facility which has maximum distance to other existing facilities with respect to the same weights and line barriers. They emphasized the rectilinear distance metric *(p=1)* in the formulation and stated that further research can be extended on different distance functions. Moreover, Amiri-Aref et al. (2013) also presented the same distance metric *(p=1)* to solve the Weber location problem with respect to a probabilistic polyhedral barrier. They transferred the barrier zone from a central point of a barrier with extreme points to a probabilistic rectangular-shaped barrier. Later, Amiri-Aref et al. (2015) utilized the "$l_1$ distance" model to deal with relocation problem with a probabilistic barrier. For solving the restricted location-relocation problem with finite time horizon, some authors developed three approaches which are well established in the literature. Dearing and Segars (2002) presented that for solving single facility location problems with barriers using rectilinear distance model ($p$=1), the result remains the same for original barriers and modified barriers.

Aneja and Parlar (1994) discussed two algorithms (Convex FR and Location with BTT) that can solve Weber Location problems with forbidden regions for $1 < p \leq 2$ when using $l_p$ distance model to calculate the distance between demand points and facility location points. Bischoff and Klamrot (2007) illustrated a solution method with respect to the problem for locating a new facility while given a set of existing facilities in the presence of convex polyhedral forbidden areas. This method at first decomposes the problem and then use genetic algorithm to find appropriate intermediate points. In the literature, only $p$=2 or Euclidean distance has been considered and it has been showed that this algorithm was very efficient. Altınel et al. (2009) considered the heuristic approach for the multi-facility location problem with probabilistic customer locations. To calculate the expected distance between probabilistic demand points and the facility locations, they used Euclidean distance, squared Euclidean distance, rectilinear distance and weighted $l_{1,2}$-norm distance models.

By virtue of network theory, we can find another way of calculating distances for location problem. Blanquero et al. (2016) represented the competitive location problem on a network with customers as the nodes in the network. They put the facilities on the edges and calculated the distances as the lengths of the shortest path from each customer to the nearest facility. Peeters (1998) discussed a new algorithm to solve the location problem by calculating distance as the sum of the lengths of arcs between two vertexes and the length between a vertex and an arbitrary point on that arc. The author used nine different models to test the new algorithms and showed much better performances than the old algorithm which calculates the distance between each node first.

Very little effort was made in the prior literature to consider concave regions thoughtfully when they estimated the distance. It may cause problems in real world application. In section 2, we will talk about the problems that could be occurred in a real-world scenario. In section 3, we will introduce our distance measurement algorithm which considers the concave regions. We will also discuss some properties as well as one observation with proof about the proposed algorithm. In section 4, we will illustrate a capacity location model which is very sensitive to the distance between regions. In Section 5, we will show how our proposed algorithm could be used to find optimal location using a test-case example. In section 6, we will consider several cases to present advantages and insights of the proposed algorithm.

## 2. Problem Description:

In the prior literature, geometric-center was widely used for distance calculation. But using the geometric-center of a region to represent itself may cause problems. Let us imagine a concave region where the geometric-center is out of the region. Now, if we use geometric-center point to represent the region then the distance between this region and another concave region might be zero. For example, in the following Figure 1, we can see Shanghai city with its outline represented as some collection concave regions. As we can see, the region for west district is a concave region; its geometric-center may fall into east district. If we transform the region a little bit, then the two geometric-centers may be the same point thus, the distance becomes zero. Therefore, in most cases, the geometric-center can represent a region well but under special situation, it will cause fatal error or inaccuracy in the final decision of a location problem. To prevent this situation, we will develop an algorithm to measure the distance between two irregular regions. We will only consider the continuous facility location problem which will allow us to set up a facility at any location. The algorithm can prevent the potential point falls out of a region itself. It also has good property within the region it represents. We will also show the relationship among the potential points we get with geometric-centers.

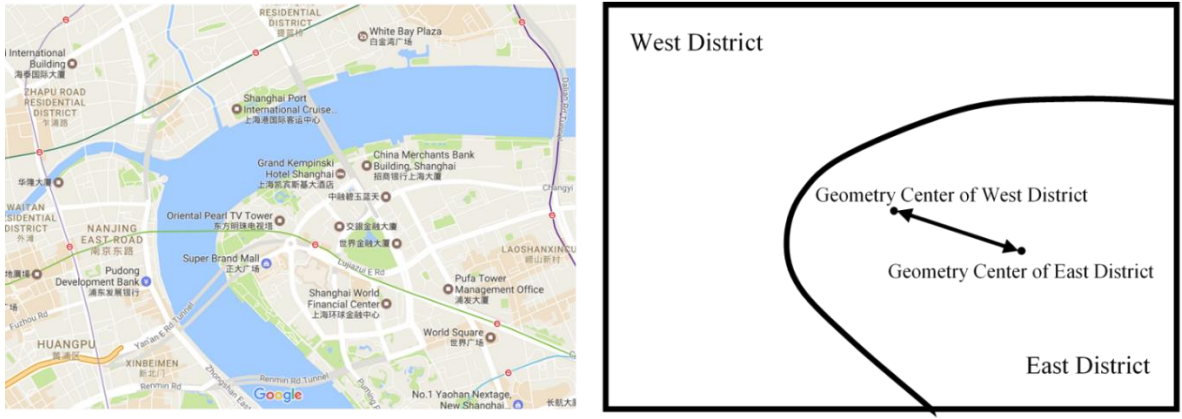

Figure 1: Map and Outline of center of Shanghai (adopted from Googlemap).

## 3  Proposed Algorithm:

3.1 Distance Measurement Algorithm:

Suppose, we have a region $S$, the point we use to represent $S$ is a point$(x, y) \in S$, such that, the value $\iint_{(x_i, y_i) \in S} \sqrt{(x - x_i)^2 + (y - y_i)^2} \, dxdy$ has a minimum value. The distance between two regions is estimated by the geometrical distance between two points which represent the regions.

It means $(x, y)$ is a point in $S$, and $(x, y)$ is the point to make the value of the function $\iint_{(x_i, y_i) \in S} \sqrt{(x - x_i)^2 + (y - y_i)^2} \, dxdy$ gains the smallest value.

We use the geometrical distance between two points we derived from the method above to estimate the distance between two regions. For example, if we have region $S_1$ and $S_2$, and we get $(x_{S_1}, y_{S_1})$ and $(x_{S_2}, y_{S_2})$ to represent $S_1$ and $S_2$, then the distance between $S_1$ and $S_2$ is $d(S_1, S_2) = \sqrt{(x_{S_1} - x_{S_2})^2 + (y_{S_1} - y_{S_2})^2}$.

3.2 Theorem:

First, it's reasonable to say that the point represents the region because $(x, y)$ is the point that makes the total travel distance needed to satisfy each demand point to gain the smallest value. The point we derived from our algorithm also equals to the geometric-center of the region if the geometric-center of the region falls into itself. We will proof the statement in the following subsections.

3.3 Observation and Proof:

If a region's geometric-center falls in the region itself, the point $(x, y)$ mentioned above is exactly the geometric-center of the region.

Since we are looking for $(x, y)$ such that, the expression $\iint \sqrt{(x - x_i)^2 + (y - y_i)^2} dxdy$ gain the minimum value, so we can remove the square root since $f(x) = \sqrt{x} (x \geq 0)$ is a strictly monotonically increasing function. Therefore, we need:

$(x, y)$ s.t. $\iint \sqrt{(x - x_i)^2 + (y - y_i)^2} dxdy$ gain the mininum value

$\Leftrightarrow (x, y)$ s.t. $\iint [(x - x_i)^2 + (y - y_i)^2] dxdy$ gain the mininum value

We can also separate the previous expression into two parts, one part only with $x$ and one part only with $y$. For $\iint (x - x_i)^2 dxdy$, since we are looking for a point $(x, y)$ such that it reaches the minimum value, we can transfer it to $\int (x - x_i)^2 dx$. So, we have:

$(x, y)$ s.t. $\iint [(x - x_i)^2 + (y - y_i)^2] dxdy$ gain the minimum value

$\Leftrightarrow (x, y)$ s.t. $\iint (x - x_i)^2 dxdy + \iint (y - y_i)^2 dxdy$ gain the minimum value

$\Leftrightarrow (x, y)$ s.t. $\int (x - x_i)^2 dx + \int (y - y_i)^2 dy$ gain the minimum value

Then, we can write our function into *Riemann sum:*

$$for\ x: \int (x - x_i)^2 dx \Leftrightarrow \lim_{n \to \infty} \sum_{i=1}^{n} (x - x_i)^2 \triangle x.$$

We obtain the mininum value of the function when $\left( \lim_{n \to \infty} \sum_{i=1}^{n} (x - x_i)^2 \triangle x \right)' = 0 \Leftrightarrow$

$\lim_{n \to \infty} \sum_{i=1}^{n} 2(x - x_i) = 0 \Leftrightarrow x = \lim_{n \to \infty} \frac{\sum_{i=1}^{n} x_i}{n}$. Since the geometric-center falls into the region, here $x$ is exactly the coordinates for geometric-center. We can get the same result for $y$. So, we proved that if a region's geometric-center falls in the region itself, the point $(x, y)$ mentioned above is exactly the geometric-center of the region.

## 4   Test Model:

In this section, we will introduce a model for testing the proposed algorithm. It's a simple capacity location problem; we have added a constraint to make it sensitive to the distance between regions.

*4.1 Parameters*

- $a_x$    demand of point $x$ at time period $a$
- $c$     fixed cost of building a facility
- $d_{xy}$   distance between points $x$ and $y$
- $L$     maximum capacity of one facility
- $M$     large cost parameter

*4.2 Variables*

- $\gamma_{xy}$   if demand at point $x$ is satisfied by facility at point $y$, then $\gamma_{xy} = 1$, else $\gamma_{xy} = 0$
- $\rho_y$   if a facility is built at $y$ at time $a$, then $\rho_y = 1$, else $\rho_y = 0$

$$Minimize \sum_y \rho_y c \quad (1)$$

$$\gamma_{xy} \leq \rho_y \quad (2)$$

$$d_{xm}\gamma_{xm} \leq d_{xy} + M(1 - \rho_y) \quad (3)$$

$$\sum_x a_x \gamma_{xy} \leq L \quad (4)$$

$$\sum_y \gamma_{xy} = 1 \quad (5)$$

The objective function (1) is a cost minimization function that aims to minimize the total cost of building facilities. Constraint (2) indicates that location $y$ can satisfy the demand at location $x$ if and only if there's a facility at location $y$. Constraint (3) restricts the demand at location $x$ to be met by its nearest facility. If there is a facility built at location $x$, then for any other point $m$, $\gamma_{xm} = 1$ if and only if $d_{xy} < d_{xm}$. This constraint will make the model very sensitive to the distance between regions. Constraint (4) shows that the total demand served by one facility can't exceed its capacity. Constraint (5) restricts that we can only set up one facility at one location.

## 5   Result Analysis:

We randomly generated a map with several concave regions shown in the following Figure 2. The proposed algorithm for distance measurement is implemented with MATLAB R 2016a. We get the geometric-center and the potential points by using our proposed algorithm. The size of the map is 200×200 pixels and we use a coordinate system on the map. We fixed the origin of coordinate at the lower left corner of the map and we also assumed that the coordinate of upper left corner is (200, 0), the coordinate of lower right corner is (0,200). Therefore, all the coordinates we assumed are integer numbers.

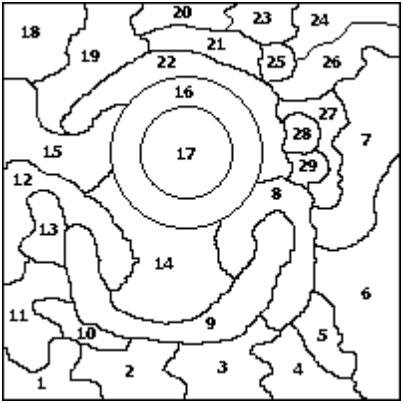

Figure 2: Test-case model map.

Table 1: Coordinates and demand of representative point for each region.

| Region Number | Geometric-center | Algorithm Center | Demand | Region Number | Geometric-center | Algorithm Center | Demand |
|---|---|---|---|---|---|---|---|
| 1 | (189,25) | (189,25) | 9 | 16 | (75,93) | (64,112) | 2 |
| 2 | (186,66) | (186,66) | 10 | 17 | (76,93) | (76,93) | 5 |
| 3 | (184,112) | (184,112) | 2 | 18 | (19,19) | (19,19) | 10 |
| 4 | (186,149) | (186,149) | 10 | 19 | (28,46) | (28,46) | 8 |
| 5 | (169,162) | (169,162) | 7 | 20 | (7,90) | (7,90) | 10 |
| 6 | (148,183) | (148,183) | 1 | 21 | (22,106) | (22,106) | 7 |
| 7 | (73,183) | (73,183) | 3 | 22 | (42,93) | (37,92) | 1 |
| 8 | (121,136) | (121,127) | 6 | 23 | (10,132) | (10,132) | 9 |
| 9 | (151,98) | (153,106) | 10 | 24 | (11,165) | (11,165) | 10 |
| 10 | (163,38) | (163,38) | 10 | 25 | (31,138) | (31,138) | 7 |
| 11 | (157,15) | (157,15) | 2 | 26 | (30,168) | (30,168) | 8 |
| 12 | (108,30) | (107,32) | 10 | 27 | (67,163) | (67,163) | 8 |
| 13 | (115,23) | (115,23) | 10 | 28 | (67,150) | (67,150) | 4 |
| 14 | (133,83) | (133,83) | 5 | 29 | (84,154) | (84,154) | 7 |
| 15 | (72,27) | (72,27) | 9 | | | | |

Table 1 shows the comparison between coordinates of representative location using geometric-center and our proposed algorithm for the point of each region. We can see that for region 8, 9, 12, 16 and 22 the representative points are different, and they are all concave regions. Now, we will use this test-case instance to find optimal location using the test case model.

The model is solved on a computer with Intel(R) Xeon(R) E3-1231 v3 CPU running at 3.4 GHz with 16 GB memory and a 64-bit operating system. For implementation, AMPL programming language is used, and CPLEX 12.6.3 solver is utilized to solve. In this instance, other parameters are assumed as follows: $c = 200$, $L = 50$, $M = 10000$. For $a_x$, values were randomly generated ranging from 1 to 10. The number in each region represents the demand at each region. Optimal facility assignments are shown in Figure 3 and Figure 4 using geometric-center and algorithm center, respectively.

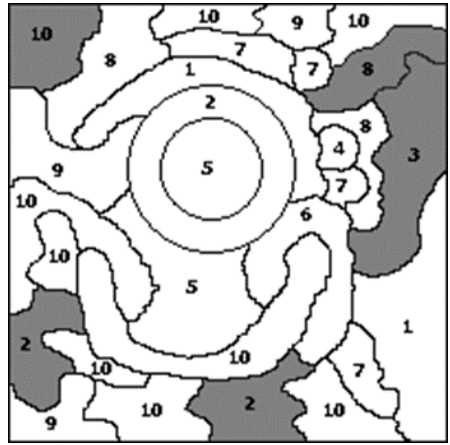

Figure 3: Optimal facility assignment using geometric-center with demand.

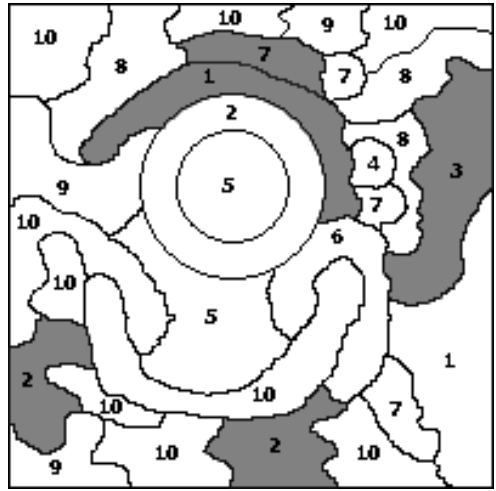

Figure 4: Optimal facility assignment using the proposed algorithm with demand.

As we can see from Table 2, two algorithms provide us the same optimal solution. However, our proposed algorithm requires less number of simplex iterations and branch-and-bound nodes. While the purpose of our algorithm is to make the distance estimation more accurate, the difference of objective function value or travel distance cannot provide better insight about it. Therefore, in the next section, we will present three test-cases to analyze the advantages of our proposed algorithm.

Comment: 1. The difference of number of simplex iterations and branch & bound nodes does not show which algorithm is better or not. It's just a coincidence. 2. It's hard to tell the difference or the advantage of our algorithm from a real case or general case.

## 6 Special Case Analysis:

In this section, we present three special cases to show further the usefulness and performance of the proposed algorithm.

Case1:
We can see from Figure 5, region 1 is a concave region and its geometric-center falls into region 2, and the actual distance between the two geometric-centers are 0. If we use this result to estimate the distances in a location problem, then there's no difference for the location of region 1 and 2. If there's another region available then the demand in that region can be either served by region 1 or region 2. Clearly, we don't want this situation to occur in real life problems. In Figure 6, we can see the representative points for region 1 falls on the boundary of region 1, which will show the difference between the region 1 and 2. If there's another region, the demand of that region can be satisfied by preference.

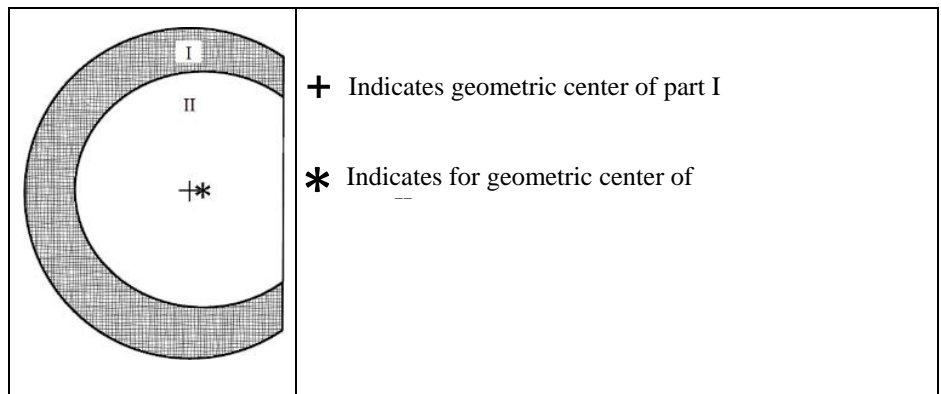

Figure 5: Geometric center for case 1.

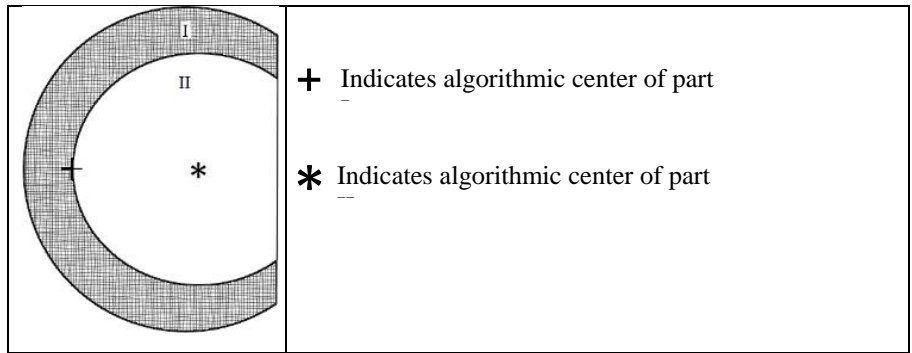

Figure 6: Algorithmic center for case 1.

Case 2:

In this case, we have 3 regions. From Figure 7, we can see that the distance between region 1 and region 2 are smaller than the distance between region 2 and region 3. This means that the demand at region 2 are more likely to be served by region 1. But in Figure 8, by our algorithm, the distance between region 2 and region 3 are smaller than the distance between region 1 and region 2. This means the demand at region 2 are more likely to be served by region 1. From the graph, we can see that it's better to let the region 2 served by region 3 since it looks closer to each other.

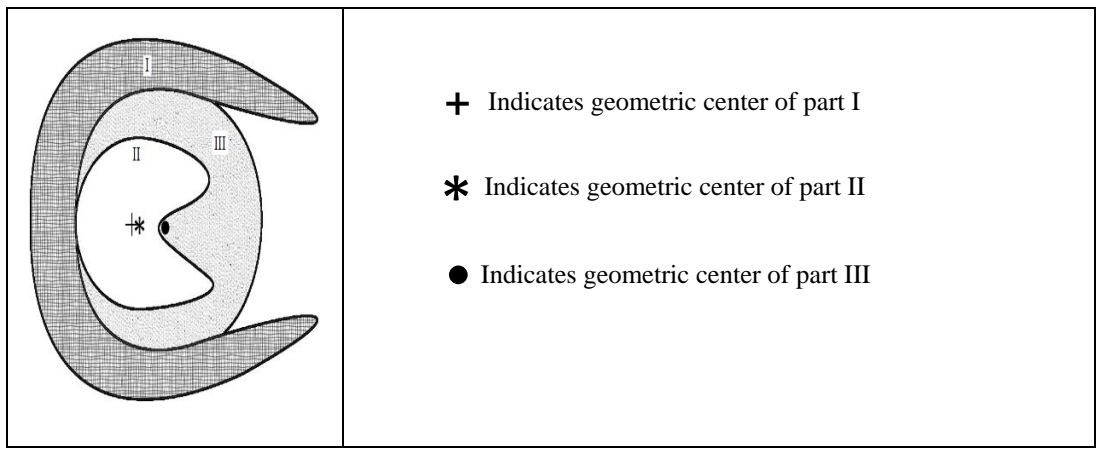

Figure 7: Geometric center for case 2.

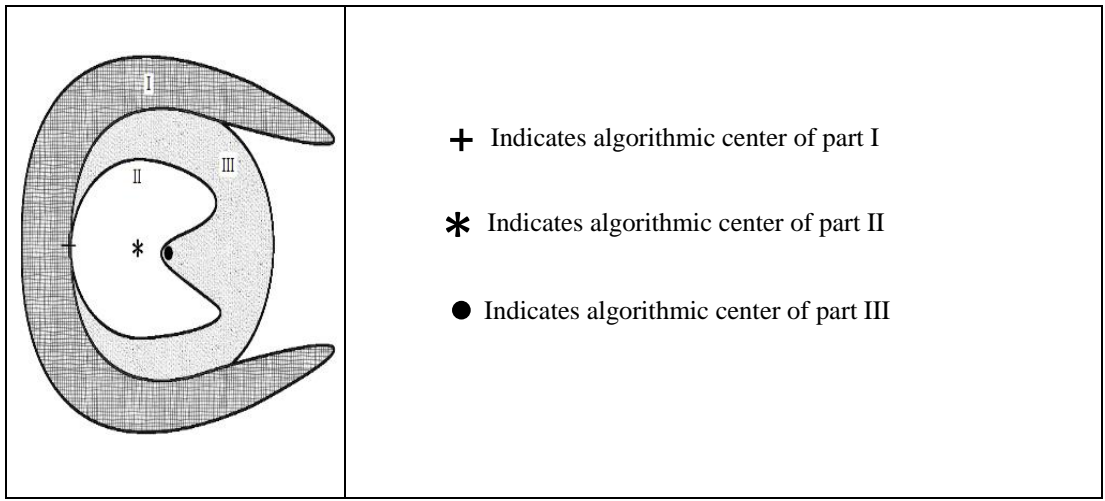

Figure 8: Algorithmic center for case 2.

Case 3:

In this case, we can see from Figure 9 and Figure 10, by both algorithms, the representative points fall together. This is a very special case for which our algorithm cannot identity the location of two regions.

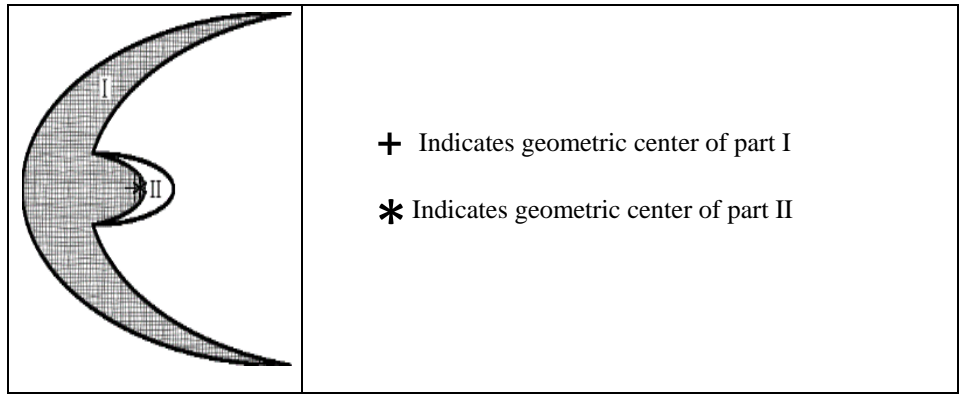

Figure 9: Geometric center for case 3.

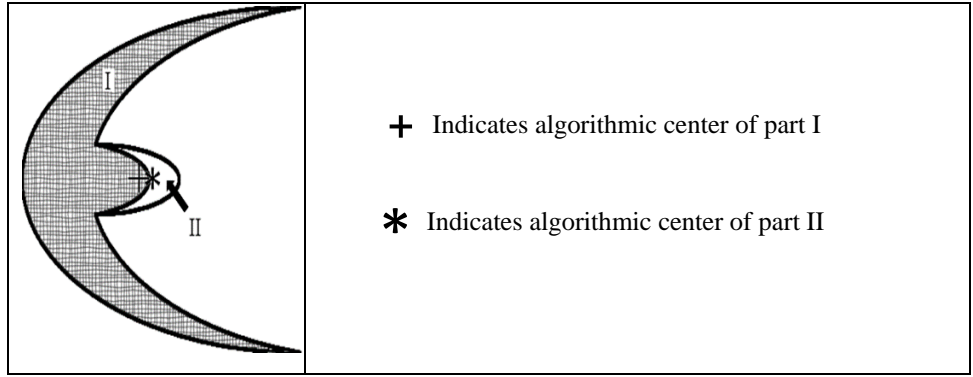

Figure 10: Algorithmic center for case 3.

## 7. Conclusion:

The proposed algorithm can identify the difference between regions when there are concave regions in the map. We have shown that under several conditions, the algorithm performs better than the geometric-center approach. When the geometric-center falls out of the region itself, our algorithm has a significant advantage to make the measurement more accurate. When the geometric-center falls into the region itself, our algorithm will keep the original feature of the geometric-center, and the algorithm center and geometric-center fall at the same point. The algorithm is also easy to implement. However, there are limitations of this algorithm; the most noticeable one is that we can only implement it on a continuous facility location problem.

**References:**

[1] Altınel, I. K., Durmaz, E., Aras, N., & Özkısacık, K. C. (2009). A location–allocation heuristic for the capacitated multi-facility Weber problem with probabilistic customer locations. *European Journal of Operational Research* 198 (3): 790-799. doi:10.1016/j.ejor.2008.10.014

[2] Amiri-Aref, M., Javadian, N., Tavakkoli-Moghaddam, R., & Baboli, A. (2013). A new mathematical model for the Weber location problem with a probabilistic polyhedral barrier. *International Journal of Production Research* 51 (20): 6110-6128. doi:10.1080/00207543.2013.796422

[3] Amiri-Aref, M., Farahani, R. Z., Javadian, N., & Klibi, W. (2015). A rectilinear distance location–relocation problem with a probabilistic restriction: mathematical modelling and

solution approaches. *International Journal of Production Research* 54 (3): 629-646. doi:10.1080/00207543.2015.1013642

[4] Amiri-Aref, M., Javadian, N., Tavakkoli-Moghaddam, R., & Aryanezhad, M. (2011). The center location problem with equal weights in the presence of a probabilistic line barrier. *International Journal of Industrial Engineering Computations* 2 (4): 793-800. doi:10.5267/j.ijiec.2011.06.002

[5] Aneja, Y. P., & Parlar, M. (1994). Algorithms for Weber Facility Location in the Presence of Forbidden Regions and/or Barriers to Travel. *Transportation Science* 28 (1): 70.

[6] Batta, R., Ghose, A., & Palekar, U. S. (1989). Locating Facilities on the Manhattan Metric with Arbitrarily Shaped Barriers and Convex Forbidden Regions. *Transportation Science 23* (1): 26-36. doi:10.1287/trsc.23.1.26

[7] Bischoff, M., & Klamroth, K. (2007). An efficient solution method for Weber problems with barriers based on genetic algorithms. *European Journal of Operational Research* 177 (1): 22-41. doi:10.1016/j.ejor.2005.10.061

[8] Blanquero, R., Carrizosa, E., G.-Tóth, B., & Nogales-Gómez, A. (2016). P-facility Huff location problem on networks. *European Journal of Operational Research 255* (1): 34-42. doi:10.1016/j.ejor.2016.04.039

[9] Dearing, P. M., & Segars Jr., R. (2002). An Equivalence Result for Single Facility Planar Location Problems with Rectilinear Distance and Barriers. *Annals of Operations Research* 111 (1-4): 89-110. doi:10.1023/A:1020945501716

[10] Francis, R. L. (1967). Letter to the Editor—Some Aspects of a Minimax Location Problem. *Operations Research* 15 (6): 1163-1169. doi:10.1287/opre.15.6.1163

[11] Love, R. F., & Dowling, P. D. (1985). Optimal Weighted lp Norm Parameters for Facilities Layout Distance Characterizations. *Management Science* 31 (2): 200-206. doi:10.1287/mnsc.31.2.200

[12] Love, R. F., & Morris, J. G. (1979). Mathematical Models of Road Travel Distances. *Management Science* 25 (2): 130-139. doi:10.1287/mnsc.25.2.130

[13] Love, R. F., & Morris, J. G. (1972). Modelling Inter-City Road Distances by Mathematical Functions. *Operational Research Quarterly (1970-1977)* 23 (1): 61. doi:10.2307/3008514

[14] Peeters, P. H. (1998). Some new algorithms for location problems on networks. *European Journal of Operational Research* 104 (2): 299-309. doi:10.1016/s0377-2217(97)00185-9